\documentstyle[12pt]{article}

\title{Perimeter and Coherence According to McCammond and Wise }
\author{Peter Brinkmann}
\date{}

\begin{document}
\maketitle
\begin{abstract}
In \cite{mw}, McCammond and Wise introduce the perimeter of
2-complexes and use it to obtain a sufficient criterion for coherence.
We present a new exposition of their main result, as well as
some applications.
\end{abstract}

\section{Introduction}
A group is called {\it coherent} if every finitely generated subgroup is
finitely presented. Coherence is a commensurability invariant, i.\ e., if
$H<G$ is a subgroup of finite index, then $G$ is coherent if and only if
$H$ is coherent.
\par

Free groups and surface groups are elementary examples
of coherent groups. Moreover, fundamental groups of 3-manifolds are coherent
(see \cite{scott}), as are mapping tori of free group automorphisms (see
\cite{fh}).
\par

The proofs of these statements heavily use specific properties of the groups
in question.  Recently, however, Wise and McCammond have developed a more
general approach to coherence that applies to a wide range of examples,
such as certain one-relator groups, some small-cancellation groups, etc.
(see \cite{mw}). The first four sections of this paper contain a new exhibition
of their sufficient criterion for coherence, and the remaining sections contain
some applications of this criterion.
\par

It is my pleasure to thank my advisor S. M. Gersten for many helpful
comments and discussions.
\par

\section{Conventions and definitions}

All complexes and maps between complexes are assumed to be combinatorial.
$f^{(n)}\in X$ denotes an $n$-cell of $X$.
Deviating from standard terminology, cells of complexes are closed,
and characteristic maps extend to the boundary of a cell. In particular,
for a 2-cell $f$ of $X$, the domain of the characteristic map
$\chi_f:D\rightarrow X$ is a polygon.
\par

For a 2-complex $X$, let $X_S$ denote its stellar subdivision. Observe that
the 2-cells of $X_S$ are triangles. For a map
$\Phi:Y\rightarrow X$ between 2-complexes, let $\Phi_\ast$ denote the induced
map at the $\pi_1$-level and let $\Phi_S:Y_S\rightarrow X_S$
denote the induced map between the stellar subdivisions of $X$ and $Y$.
For an edge $e$ of $X$, the {\it star} $St(e)$ of $e$ is the collection
of triangles of $X_S$ that are adjacent to the edge of $X_S$
corresponding to $e$.
\par

\subsection{Definition}
A {\it weight function} on a 2-complex $X$ is a map $w:\{$triangles of $X_S\}
\rightarrow {\bf R}$. The weight $w(f)$ of a 2-cell $f$ is the sum of
the weights of the triangles in its stellar subdivision.
The {\it standard weight function}
assigns the weight $1$ to every triangle in $X_S$.
Unless stated otherwise, weight functions in this paper will take values in
the nonnegative integers. 
\par
Given a weight function on $X$ and a map $\Phi:Y\rightarrow X$, the
{\it missing weight} of an edge $e$ of $Y$ is the sum of the weights of those
triangles in the star of $\Phi(e)$ that are not contained in the image of
the star of $e$. The {\it missing weight}
of $\Phi$ is the sum of the missing weights
of the edges of $Y$. In formulas, we have the double sum
$$M(\Phi)=\sum_{e^{(1)}\in Y}\;
\sum_{f^{(2)}\in St(\Phi_S(e))\setminus \Phi_S(St(e))} w(f).$$
The missing weight of a subcomplex $Y\subset X$ is the missing weight of
the inclusion map. 
\par
Remark: We use the term ``missing weight'' instead of ``perimeter'' since
the latter has caused some confusion.
\par

Example: If $X$ is a combinatorial 2-manifold without boundary (equipped
with the standard weight function), then
the missing weight of each 1-cell of $X$ equals 2.
If $X$ is the 2-skeleton of the usual cubulation of ${\bf R}^3$,
then the missing weight of each 1-cell is 4, and the missing weight of 
each 2-cell is 12.
\par

A 2-cell $f$ of $X$ with characteristic map $\chi_f:D\rightarrow X$
is said to have {\it exponent} $n$ if there is a closed
path $w$ in $X^{(1)}$ such that $\chi_f\mid_{\partial D}$ spells the word
$w^n$ and if $n$ is maximal. The polygon $D$ is symmetric under rotations by
${2\pi\over n}$.\par
The {\it packet} $\tilde f$ of $f$ is the complex obtained by taking $n$
copies of $D$ and gluing them along the boundary. There is a well defined map
$\chi_{\tilde f}:\tilde f\rightarrow X$ whose restriction to the $i$th copy
of $D$ is $\chi_f$ precomposed with a rotation by ${2\pi\over n}$.
\par
A map $\Phi:Y\rightarrow X$ is said to be {\it packed} if for every 2-cell
$f$ of $X$ and every lift $\gamma$ of $\chi_f$ to $Y$ there exists a 
lift of $\chi_{\tilde f}$ to $Y$ extending $\gamma$.
This condition is vacuous for any 2-cell whose characteristic map does not
factor through $Y$. By attaching 2-cells to $Y$ and extending $\Phi$
appropriately, we can assume that $\Phi$ is packed.
Note that for
nonnegative weights this does not increase the missing weight and that it does
not change anything at the $\pi_1$-level.
\par

\subsection{Example}\label{perex}
Let $X$ be the presentation complex of $<x\mid x^n>$. $X$ has exactly one
2-cell with characteristic map $\chi_f:D\rightarrow X$. Let $\Phi_1:X^{(1)}
\rightarrow X$ denote the inclusion of the 1-skeleton, and let $\Phi_2$
be the restriction of $\chi_f$ to the boundary of $D$.
\par
Then the missing weight of $\Phi_1$ is exactly the weight of $f$, and packing
$\Phi_1$ yields a map of missing weight $0$. The missing weight of $\Phi_2$ is
$n\cdot w(f)$, extending $\Phi_2$ to the interior of $D$ reduces the
missing weight by $w(f)$, and packing $\Phi_2$ yields a map of missing weight
$0$ since each copy of $D$ reduces the missing weight by $w(f)$.
\par

\section{Reductions}

Roughly speaking, the missing weight
 of a map $\Phi:Y\rightarrow X$ will play the
role of a complexity function measuring to what extent $\Phi$ fails to
be $\pi_1$-injective.
\par
In order to show that a finitely generated subgroup
$H < \pi_1X$ is finitely presented, we will start with some finite complex
$Y$ and a map $\Phi:Y\rightarrow X$ such that $\Phi_\ast(\pi_1Y)=H$. If
$\Phi_\ast$ has nontrivial kernel, we will attach a suitable 2-cell to $Y$
in order to make the kernel smaller. If we can do this in such a way that
the missing weight goes down with each reduction of the kernel, this process
terminates after finitely many steps, proving that $H$ is finitely presented.
If this works for any finitely generated $H<\pi_1X$, it follows that $\pi_1X$
is coherent, and we will give sufficient conditions for this to happen.

\subsection{Definition}\label{reddef}
Fix a packed map $\Phi:Y\rightarrow X$ and a 2-cell $f$ of $X$ with
characteristic map $\chi_f:D\rightarrow X$.
Let $P$ be an interval subdivided into no more than $|\partial D|$ edges,
where $|\partial D|$ denotes the number of faces of $D$.
\par
A pair $(\rho:P\rightarrow D,\rho':P\rightarrow Y)$ of immersed paths
is called a {\it reduction of $\Phi$} if the
following conditions are satisfied:
\begin{enumerate}
\item $\rho$ and $\rho'$ fit into the following commutative diagram:
\begin{eqnarray*}
P & \stackrel{\rho'}{\longrightarrow} & Y \\
\rho\downarrow &\nearrow\!\!\!\!\!\!{\scriptstyle \setminus} &
\downarrow\Phi \\
D & \stackrel{\chi_f}{\longrightarrow} & X\\
\end{eqnarray*}
\item There is no map $D\rightarrow Y$ that will fit into the diagram above.
\end{enumerate}

If $|P|<|\partial D|$, the reduction $(\rho, \rho')$ is said to be
{\it incomplete}, in which case we define another immersed path
$\sigma:S\rightarrow \partial D$ such that $S$ is an interval subdivided into
$|S|=|\partial D|-|P|$ edges and $\rho(P)\cup \sigma(S)=\partial D$.
$\sigma$ is called a {\it complement} of $\rho$, and it is unique up to
orientation.
\par

If $|P|=|\partial D|$, the reduction $(\rho,\rho')$ is called {\it complete}.
Finally, $(\rho,\rho')$ is said to be {\it maximal} if there is no reduction
$(\tau,\tau')$ with the property that $\rho$ and $\rho'$ are proper
subpaths of $\tau$ and $\tau'$, respectively.

\par

Roughly speaking, reductions allow us to attach packets of 2-cells to $Y$, thus
(under suitable hypotheses)
reducing the missing weight of $\Phi$. More precisely, we have the
following
\subsection{Construction}\label{constr}
Fix a map $\Phi:Y\rightarrow X$ and a 2-cell $f$ of $X$ with characteristic
map $\chi_f:D\rightarrow X$. If the pair
$(\rho:P\rightarrow D,\rho':P\rightarrow Y)$ is a reduction, then we construct
a new map $\Phi^+:Y^+\rightarrow X$ in the following way:\par
If $(\rho,\rho')$ is a complete reduction, then the endpoints of $\rho$ are
necessarily the same, even though the endpoints of $\rho'$ in $Y$ may not be
equal.
In this case, identify the endpoints of $\rho'$ in $Y$, obtaining a new complex
$Y'$. $\Phi$ factors through a map $\Phi':Y'\rightarrow X$. Abusing notation,
we still refer to the path $P\rightarrow Y \rightarrow Y'$ as $\rho'$.
In the other case, i.\ e.,
if the endpoints of $\rho'$ are already equal or if the
reduction is incomplete, let $Y'=Y$ and $\Phi'=\Phi$. In any case, we have
$M(\Phi)=M(\Phi')$.\par
Next, we define a complex $Y''$ by amalgamating $D$ and $Y'$
along $P$, i.\ e., $Y''=Y'\cup_P D=(Y'\coprod D)/\!\sim$,
where two elements $y\in Y'$
and $x\in D$ are equivalent precisely if there exists some $t\in P$ such that
$\rho(t)=x$ and $\rho'(t)=y$. We can think of $Y'$ and $D$ as being
subcomplexes of $Y''$, and we define a map $\Phi'':Y''\rightarrow X$ by
$\Phi''|_{Y'}=\Phi'$ and $\Phi''|_D=\chi_f$.\par
Finally, we obtain $\Phi^+:Y^+\rightarrow X$ by packing
$\Phi'':Y''\rightarrow X$. Note that $\Phi^+_\ast(\pi_1Y^+)=\Phi_\ast(\pi_1Y)$.
\par

\subsection{Lemma} \label{incompred}
{\it Let $\Phi:Y\rightarrow X$ be a packed map, and let $f$ be a 2-cell of $X$
with characteristic map $\chi_f:D\rightarrow X$. If $(\rho:P\rightarrow D,
\rho':P\rightarrow Y)$ is a maximal incomplete reduction, 
then $$ M(\Phi^+)=M(\Phi)+M(\chi_f\circ\sigma)-n\cdot w(f),$$
where $n$ is the exponent of $f$, $\sigma$ is a complement of $\rho$ and
$\Phi^+$ is the map constructed in \ref{constr}.}
\par
Proof: Choose a word $w$ in $X^{(1)}$ such that $\chi_f|_{\partial D}$
spells the word
$w^n$, where $n$ is the exponent of $f$. If $\rho':P\rightarrow Y$
has a closed loop as a subpath, say $\rho''$, then $\Phi\circ\rho''$
necessarily spells a
power of a conjugate of $w$ in $X^{(1)}$ since $\Phi:Y\rightarrow X$ is
combinatorial. This implies that the reduction $\rho$ can be extended, which
contradicts the maximality assumption. Hence, $\rho'$ contains no closed loop.
\par
Since $\rho'$ contains no closed loop, the packet $\tilde f$
is a subcomplex of $Y^+$, and
we have $M(\Phi^+|_{Y^+\setminus \{\rm interior\ of\ 2-cells\ of\ \tilde f\}})=
M(\Phi)+M(\chi_f\circ\sigma)$. Since there is no map $D\rightarrow Y$ that
will fit into the diagram in definition \ref{reddef},
gluing back a 2-cell reduces the missing weight
exactly by $w(f)$ (cf. example \ref{perex}). $\Box$

\subsection{Lemma} \label{compred}
{\it Let $\Phi:Y\rightarrow X$ be a packed map, and let $f$ be a 2-cell of $X$
with characteristic map $\chi_f:D\rightarrow X$. If $(\rho:P\rightarrow D,
\rho':P\rightarrow Y)$ is a complete reduction, then
$$M(\Phi^+)\le M(\Phi)-w(f),$$
where $\Phi^+$ is the map constructed in \ref{constr}.}
\par

Proof: Since $\rho$ is complete, all the 1-cells in $Y^+$ correspond to
1-cells in $Y$, which implies that $M(\Phi)=
M(\Phi^+|_{Y^+\setminus \{\rm interior\ of\ 2-cells\ of\ \tilde f\}})$.
Since there is no map $D\rightarrow Y$ that will fit into the commutative
diagram in definition \ref{reddef}, gluing the first
2-cell of $\tilde f$ back reduces the missing weight by $w(f)$,
which proves the claim. Note that (as opposed to the previous proof) $\rho'$
may contain closed loops, in which case gluing back subsequent 2-cells may
not result in a further reduction of missing weight.
$\Box$\par

\section{A sufficient condition for coherence}
Let $\alpha:Q\rightarrow X$ be a contractible immersed loop in $X^{(1)}$. 
Then there exists a sequence of paths $\alpha_0=\alpha,\cdots,
\alpha_k=const$ such that $\alpha_{i+1}$ is obtained from $\alpha_i$
either by tightening, i.\ e., removing all subpaths of the form $e\bar{e}$, or by
a homotopy across a 2-cell in the following way: A homotopy across a 2-cell $f$
with characteristic map $\chi_f:D\rightarrow X$ uses a maximal reduction
$(\rho:P\rightarrow D,\rho':\rightarrow Q)$ such that
\begin{eqnarray*}
P & \stackrel{\rho'}{\longrightarrow} & Q \\
\rho  \downarrow & & \downarrow  \alpha_i \\
D & \stackrel{\chi_f}{\longrightarrow} & X\\
\end{eqnarray*}
commutes. Note that we can think of $P$ as being a
subcomplex of $Q$. If $\rho$ is
complete, then $\alpha_{i+1}$ is constructed from
$\alpha_i$ by removing $P$ from $Q$. If $\rho$ is
incomplete, then we replace the subpath $\chi_f\circ\rho$ with its complement
$\chi_f\circ\sigma$ (notation as in definition \ref{reddef}; we have to
choose the correct the orientation of $\sigma$ in order to obtain a
continuous path).

\subsection{Definition}\label{pathred}
The complex $X$ is said to have the {\it path reduction property} if for each
contractible immersed loop $\alpha$ there is a sequence $\alpha_0=\alpha,
\cdots,\alpha_k=const$ as above such that for every incomplete reduction
$(\rho,\rho')$ involving a 2-cell $f$ we have the inequality 
$$M(\chi_f\circ\sigma)\le n\cdot w(f),$$
where $n$ is the exponent of $f$ and $\sigma$ is a complement of $\rho$ (cf.
previous paragraph).
\par

\subsection{Theorem}\label{main}
{\it Let $X$ be a 2-complex enjoying the path reduction property for some weight
function $w$. If the missing weight (with respect to $w$) of each edge of
$X$ is finite and if the weight of each 2-cell of $X$ is strictly positive,
then $\pi_1X$ is coherent.}
\par
\bigskip

Proof: Let $H<\pi_1X$ be a finitely generated subgroup. We want to show that
$H$ is finitely presented. There exists
a finite 1-complex $Y$ and a map $\Phi:Y\rightarrow X$ such that 
$\Phi_\ast(\pi_1Y)=H$ and $\Phi$ is an immersion on the 1-skeleton (see
\cite{tfg} for details).
\par
If $\Phi$ is not $\pi_1$-injective, then there exists an essential loop
$\beta$ in $Y$ such that $\alpha=\Phi\circ\beta$ is contractible in $X$.
Choose a sequence $\alpha_0=\alpha,\cdots,\alpha_k=const$ as above,
using the path reduction property of $X$. The idea is to add finitely many
2-cells to $Y$ in such a way that the loops $\alpha_0,\cdots,\alpha_k$
and the homotopies between them lift to the new complex.
\par
Assume inductively that the loops $\alpha_0,\cdots,\alpha_i$ and the
homotopies between them have already been lifted to $\beta_0=\beta,\cdots,
\beta_i$ in $Y$.
If $\alpha_{i+1}$ is obtained from $\alpha_i$ by tightening, then
$\alpha_{i+1}$ lifts to $Y$ since $\Phi$ is an immersion on the 1-skeleton.
\par
If a homotopy across a 2-cell fails to lift to $Y$, then the corresponding
reduction $(\rho:P\rightarrow D,\rho':P\rightarrow Q)$
of $\alpha_i:Q\rightarrow X$ gives rise to a reduction of $\Phi$ in the
following way: The path $\alpha_i$ in the commutative diagram
\begin{eqnarray*}
P & \stackrel{\rho'}{\longrightarrow} & Q \\
\rho  \downarrow & & \downarrow  \alpha_i \\
D & \stackrel{\chi_f}{\longrightarrow} & X\\
\end{eqnarray*}
factors through $Y$ by our inductive hypothesis, inducing a
reduction $(\rho,\beta_i\circ\rho')$ of $\Phi$:
\begin{eqnarray*}
  P & \stackrel{\rho'}{\longrightarrow} Q
\stackrel{\beta_i}{\longrightarrow} & Y  \\
\rho  \downarrow & & \downarrow  \Phi \\
D & \stackrel{\chi_f}{\longrightarrow} & X\\
\end{eqnarray*}
Extend this to a maximal reduction of $\Phi$ and form the map
$\Phi^+:Y^+\rightarrow X$ as in \ref{constr}.
Now the loop $\alpha_{i+1}$ lifts to $Y^+$, as does the homotopy taking
$\alpha_i$ to $\alpha_{i+1}$.
\par
If the reduction is incomplete, then lemma \ref{incompred} implies that
$M(\Phi^+)=M(\Phi)+M(\chi_f\circ\sigma)-n\cdot w(f)$, and the path reduction
property guarantees that $M(\Phi^+)\le M(\Phi)$. If the reduction is complete,
then, by lemma \ref{compred},
$M(\Phi^+)\le M(\Phi)-w(f) < M(\Phi)$ since $w(f)>0$ by assumption.
$\Phi^+$ may not be an immersion on the 1-skeleton, but we can correct this by
folding edges (see \cite{tfg}) without increasing the missing weight, which
completes the inductive step.
\par
Observe that we cannot arrive at a constant path unless at least
one of the reductions of $\Phi$ is complete, so we end up with a finite
complex $Y'$ and a map $\Phi':Y'\rightarrow X$ such that
$\Phi'_\ast(\pi_1Y')=H$, $\Phi'$ is an immersion on the 1-skeleton,
and $M(\Phi')<M(\Phi)$.
\par
Since we are using nonnegative integer weights, we can only repeat this
process finitely many times before we arrive
at a $\pi_1$-injective map $\Phi'':Y''\rightarrow X$ with finite domain
$Y''$, which implies that $H$ is finitely presented. Since $H$ was arbitrary,
this implies that $\pi_1X$ is coherent. $\Box$

\section{Applications}
The applications listed in this section can be found in \cite{mw}; we give a
unified approach to them using matchings in bipartite graphs.
We will use the following concepts and results graph theory (see
\cite{graphs}, sec.\ 3.1 for details):
\subsection{Matchings in bipartite graphs}
All graphs are assumed to be finite.
A graph $G$ is called {\it bipartite} if its vertex set $V$ can be expressed
as the disjoint union of two nonempty sets $V_1, V_2$ such that every edge
of $G$ connects an element of $V_1$ to an element of $V_2$. In this case,
a {\it matching of $V_1$ into $V_2$} is a set $M$ of pairwise disjoint edges
of $G$ such that every vertex in $V_1$ is contained in some edge in $M$.
\par
A bipartite graph $G$ with bipartition $V_1, V_2$ is said to have the
{\it matching property with respect to $V_1$} if for any subset $S\subset
V_1$ the number of vertices in $V_2$ 
adjacent to $S$ is at least the number of elements of $S$.
The following theorem holds (see \cite{hall,graphs}):
\subsection{Hall's matching condition}\label{match}
{\it If $G$ is a bipartite graph with bipartition $V_1, V_2$, then $G$ has a
matching of $V_1$ into $V_2$ if and only if $G$ has the matching property with
respect to $V_1$.}
$\Box$

\bigskip
We now generalize Hall's theorem.
Given a function $m$ that assigns a positive integer to every vertex in $V_1$,
we call a set $M$ of edges an {\it m-matching} if each
element of $V_2$ is contained in at most one element of $M$, and if for
every vertex $v$ in $V_1$, $M$ contains exactly $m(v)$ edges emanating from
$v$. We call $m(v)$ the {\it multiplicity} of $v$.
\par
The graph $G$ is said to have the {\it m-matching property }
if for every subset $W$ of $V_1$, the number of elements of $V_2$
adjacent to elements of $W$ is at least as large as the sum of the
multiplicities of the elements of $W$.
\par

\subsection{The $m$-matching theorem}\label{mmatch}
{\it Let $G$ be a bipartite graph with bipartition $V_1, V_2$, and let $m$
be a multiplicity function on $V_1$. Then $G$ allows an $m$-matching if and
only if $G$ has the $m$-matching property.}
\par
Proof: Let $G'$ be a graph whose vertex set is the disjoint union of $V_2$
and $V_1'$, where $V_1'$ contains vertices $v^x_1\cdots v^x_{m(x)}$ for
every $x\in V_1$. For every edge $(x,v)$ connecting $x\in V_1$ and $v\in V_2$
we choose edges $(v^x_1,v)\cdots (v^x_{m(x)},v)$.\par
Now $G'$ has the matching property if and only if $G$ has the $m$-matching
property, and $G'$ has a matching of $V_1'$ into $V_2$ if and only if
$G$ has an $m$-matching, so the claim follows from Hall's matching
theorem. $\Box$

\section{Multiplicities, matchings and coherence}
Let ${\cal P}=<X\ |\ R>$ be a finite presentation. The {\it
incidence graph} of $\cal P$ is the bipartite graph $G$ with vertex set
$V=R\sqcup X$ such that $r\in R$ and $x\in X$ are connected by an edge if
and only if $x$ occurs in the spelling of $r$ (relators are assumed to be
cyclically reduced). Given some multiplicity function $m$ on $R$, $\cal P$
is said to have the $m$-matching property if $G$ has the $m$-matching
property.
\par
The following theorem is the main result of this section. It shows that,
for the right choice of multiplicities, the $m$-matching property implies
coherence. This theorem unifies the ideas behind several theorems in
\cite{mw}.
\par

\subsection{Theorem}\label{mainapp}
{\it Let the group $G$ be given by the finite presentation $\cal P$
satisfying the $m$-matching property for some multiplicity function $m$.
If we can reduce any word representing $1$ to the empty word by cyclic
reduction or by replacing a subword of a relator $r$ by its complement
$\sigma$ in such a way that
\begin{equation}
n\cdot m(r)\geq |\sigma|,\label{appineq}
\end{equation}
then $G$ is coherent. Here $n$ is the exponent of $r$.}
\par \bigskip
Proof: Let $G$ be the incidence graph of ${\cal P}=<X\ |\ R>$, and let $Y$
be the presentation complex of $\cal P$.
By hypothesis, $G$ has an $m$-matching $M$.
Pick some edge in $M$, connecting some $r\in R$ to some $x\in X$.
Now consider the 2-cell $f$ of $Y$ corresponding to $r$. At
least one of the triangles in the stellar subdivision of $f$ is adjacent to
the 1-cell corresponding to $x$. Assign the weight one to one such triangle.
Repeat this for all edges in $M$, and assign the weight zero to all remaining
triangles.
\par
Since each element of $X$ belongs to at most one edge in $M$, the missing weight
of each edge is at most one, so we have $|\sigma|\geq M(\sigma)$ for any
path $\sigma$. Moreover, the weight of each 2-cell is exactly the multiplicity
of the corresponding relator. Since inequality (\ref{appineq}) holds,
$Y$ has the path
reduction property. Since the weights of all 2-cells are positive, this shows
that $G$ is coherent. $\Box$
\par
\bigskip

When applying this theorem to a presentation $\cal P$, one typically begins
with a multiplicity function $m$ that satisfies (\ref{appineq}), then one
checks whether $\cal P$ has the $m$-matching property. Hence it is
advantageous to choose $m$ as small as possible.
We list some reasonable choices for certain classes of presentations.

\subsubsection{Dehn presentations}
A presentation $\cal P$ is said to be a {\it Dehn presentation} if Dehn's
algorithm solves the word problem for $\cal P$ (Dehn's algorithm solves the
word problem for $\cal P$ if any cyclically reduced word representing the
identity element contains a subword $u$ of a relator $r$ such that
$|u|>{{|r|}\over 2}$).
\par
For a Dehn presentation $\cal P$ and a relator $r$ of $\cal P$,
let $n$ be the exponent of $r$ and let $w$ be a
cyclically reduced word $w$ such that $r=w^n$. We define
$m(r) =\lfloor {{|w|-1}\over 2} \rfloor$ if $n=1$,
$m(r) =\lfloor {{|w|}\over 2} \rfloor$ if $n=2$, and
$m(r) =\lfloor {{|w|+1}\over 2} \rfloor$ if $n\geq 3$.
\par

This choice of multiplicities satisfies (\ref{appineq}) because it implies
$n\cdot m(r)\geq \lfloor {{|r|-1}\over 2} \rfloor$.
Hence, by theorem \ref{mainapp}, a group given by a Dehn presentation
is coherent if it has the matching property with respect to this choice of
multiplicities.

\subsubsection{Remark}
Since presentations of type $C(4)-T(5)$ and $C(3)-T(7)$ define hyperbolic
groups, it is natural to ask whether such presentations are Dehn presentations.
The following two examples show that in general, the answer is no: Let
${\cal P}_1 = <a, b, c, d, r, s \ |\ ab\bar{a}r, \bar{r}\bar{b}cs,
\bar{s}d\bar{c}\bar{d}>$ and ${\cal P}_2 = <a, b, c, d, r, s, t, u, v \ |\ 
abr,\- \bar{r}\bar{a}s,\- \bar{s}\bar{b}t,\- \bar{t}cu,\- \bar{u}dv,\-
\bar{v}\bar{c}\bar{d} >$.
${\cal P}_1$ and ${\cal P}_2$ are of type $C(4)-T(5)$ and $C(3)-T(7)$,
respectively, and both ${\cal P}_1$ and ${\cal P}_2$ are presentations of the
fundamental group of the closed surface of genus two, so their star graph
has exactly one cycle. This implies that  -- up to inversion and conjugacy --
there is exactly one shortest boundary label $w$ of van Kampen diagrams with
exactly one interior vertex. It is easy to check that $w$ cannot be shortened
by means of a homotopy across a relator disc, which shows that
${\cal P}_1$ and ${\cal P}_2$ are not Dehn presentations.
\par

\subsubsection{Small cancellation conditions}
A presentation is said to have property $P$ if every piece has length one and
no relator is a proper power (see \cite{gs}).\par
Since a nonempty cyclically reduced word representing 1 in a
$C(6)$ presentation contains a complement of no more than three pieces, 
inequality (\ref{appineq}) holds if we assign the multiplicity three to all the
relators in a $C(6)-P$ presentation. Hence we recover a special case of
theorem 9.4 of \cite{mw}.
\par
A nonempty cyclically reduced word representing 1 in a $C(4)-T(4)$ presentation
contains a complement of no more than two pieces, so inequality (\ref{appineq})
holds if we assign the multiplicity two to all the relators in a $C(4)-T(4)-P$
presentation. In this case we recover a special case of theorem 9.5 of
\cite{mw}.
\par

\subsubsection{$\lambda$-presentations}
A presentation is said to be a {\it $\lambda$-presentation} if every
nontrivial word representing the neutral element contains a subword of some
relator $r$ of length strictly greater than $(1-\lambda)|r|$. For example,
Dehn presentations are ${1\over 2}$-presentations by definition, and B.B.
Newman's spelling theorem shows that presentations of the form $<X\ |\ w^n>$
are ${1\over n}$-presentations if $n\geq 2$.
\par
For a relator $r$ with exponent $n$ in a $\lambda$-presentation, let $k$ be
the {\it largest} integer satisfying $k<\lambda |r|$. Then we choose $m(r)$
to be the {\it smallest} positive integer satisfying $n\cdot m(r)\geq k$.
By construction, inequality (\ref{appineq}) holds.
\par
For example, if ${\cal P}=<X\ |\ w^n>$ for some $n\geq 2$, we have
$k=|w|-1$, and $m(r)\geq {{|w|-1}\over n}$.
The incidence graph of $\cal P$ clearly has the $m$-matching property if
the number of generators occuring in $w$ is at least ${{|w|-1}\over n}$.
\par In particular,
the group given by $\cal P$ is coherent if $n\geq |w|-1$, so we
recover theorem 7.3 in \cite{mw}.
\par
\medskip

\bibliographystyle{alpha}
\bibliography{per}
\bigskip
{\sc Department of Mathematics, University of Utah\\
Salt Lake City, UT 84112, USA\\}
\\
{\it E-mail:} brinkman\@@math.utah.edu

\end{document}